\documentclass[12pt]{article}
\usepackage{mathrsfs}
\usepackage{epic,eepic,epsf,epsfig}
\usepackage{amsfonts,srcltx,mathrsfs}
\usepackage{multirow}
\usepackage{amsmath}
\usepackage{amssymb}
\usepackage{amsbsy}
\usepackage{graphicx}
\usepackage{amsfonts}
\usepackage{color}
\newtheorem{lem}{Lemma}[section]%
\newtheorem{theorem}[lem]{Theorem}%
\newtheorem{cor}[lem]{Corollary}%
\newtheorem{prop}[lem]{Proposition}%

\def\a{\alpha}

 \def\O{\Omega} \def\G{\Gamma}

\def\di{\bigm|}  
\def\nd{\mathrel{\bigm|\kern-.7em/}}

\def\f{\noindent}

\def\PSL{\hbox{\rm PSL}}\def\PSU{\hbox{\rm PSU}}
  \def\GL{\hbox{\rm GL}} 
\def\PSp{\hbox{\rm PSp}}\def\P\GammaL{\hbox{\rm P\Gamma L}}

\def\Out{\hbox{\rm Out}}
\def\Aut{\hbox{\rm Aut}}

\def\soc{\hbox{\rm soc}}

\def\rad{\hbox{\rm rad}}

\newcommand{\qed}{\mbox{\raisebox{0.7ex}{\fbox{}}} \vspace{4truemm}}
\def\mz{{\mathbb Z}}

\begin{document}
\title{Pentavalent symmetric graphs admitting transitive non-abelian characteristically simple groups}

\author{Jia-Li Du, Yan-Quan Feng\footnotemark\\
{\small\em Department of Mathematics, Beijing
Jiaotong University, Beijing 100044, China}}

\footnotetext[1]{Corresponding author. E-mails:
JiaLiDu@bjtu.edu.cn, yqfeng$@$bjtu.edu.cn}

\date{}
 \maketitle

\begin{abstract}

Let $\G$ be a graph and let $G$ be a group of automorphisms of $\G$. The graph $\G$ is called {\it $G$-normal} if $G$ is normal in the automorphism group of $\G$. Let $T$ be a finite non-abelian simple group and let $G = T^l$ with $l\geq 1$. In this paper we prove that if every connected pentavalent symmetric $T$-vertex-transitive graph is $T$-normal, then every connected pentavalent symmetric  $G$-vertex-transitive graph is $G$-normal. This result, among others, implies that every connected pentavalent symmetric $G$-vertex-transitive graph is $G$-normal except $T$ is one of $57$ simple groups. Furthermore, every connected pentavalent symmetric $G$-regular graph is $G$-normal except $T$ is one of $20$ simple groups, and every connected pentavalent $G$-symmetric graph is $G$-normal except $T$ is one of $17$ simple groups.

\bigskip
\f {\bf Keywords:} Vertex-transitive graph, symmetric graph, Cayley graph, regular permutation group, simple group.
\end{abstract}

\section{Introduction}

Throughout this paper, all groups and graphs are finite, and all graphs are simple and undirected. Denote by $\mz_n$, $D_n$, $A_n$ and $S_n$ the cyclic group of order $n$, the dihedral group of order $2n$, the alternating group and the symmetric group of degree $n$, respectively. Let $G$ be a permutation group on a set $\O$ and let $\a\in \O$. Denote by $G_{\a}$ the stabilizer of $\a$ in $G$, that is, the subgroup of $G$ fixing the point $\a$. We say that $G$ is {\em semiregular} on $\O$ if $G_\a=1$ for every $\a \in \O$, and {\em regular} if it is semiregular and transitive. For a graph $\G$, we denote its vertex set and automorphism group by $V(\G)$ and $\Aut(\G)$, respectively. The graph $\G$ is said to be {\em $G$-vertex-transitive} or {\em $G$-regular} for $G\leq \Aut(\Gamma)$ if $G$ acts transitively or regularly on $V(\G)$ respectively, and {\em $G$-symmetric} if $G$ acts transitively on the arc set of $\G$ (an arc is an ordered pair of adjacent vertices). In particular, $\G$ is {\em vertex-transitive} or {\em symmetric} if it is $\Aut(\G)$-vertex-transitive or $\Aut(\G)$-symmetric, respectively. A graph $\G$ is said to be {\em $G$-normal} for $G\leq \Aut(\G)$ if $G$ is normal in $\Aut(\G)$.

For a non-abelian simple group $T$, $T$-vertex-transitive graphs have received wide attentions, specially for the two extreme cases: $T$-symmetric graphs and $T$-regular graphs. It was shown in~\cite{DFZh} that a connected pentavalent symmetric $T$-vertex-transitive graph $\G$ is either $T$-normal or $\Aut(\G)$ contains a non-abelian simple normal subgroup $L$ such that $T\leq L$ and $(T,L)$ is one of $58$ possible pairs of non-abelian simple groups.

A $T$-regular graph is also called a {\em Cayley graph} over $T$, and the Cayley graph is called {\em normal} if it is $T$-normal. Investigation of Cayley graphs over a non-abelian simple group is currently a hot topic in algebraic graph theory. One of the most remarkable achievements is the complete classification of connected trivalent symmetric non-normal Cayley graphs over non-abelian simple groups. This work was began in 1996 by Li~\cite{CHLi}, and he proved that a connected trivalent symmetric Cayley graph $\G$ over a non-abelian simple group $T$ is either normal or $T=A_5$, $A_7$, $\PSL(2,11)$, $M_{11}$, $A_{11}$, $A_{15}$, $M_{23}$, $A_{23}$ or $A_{47}$. In 2005, Xu {\em et al}~\cite{XFWX2005} proved that either $\G$ is normal or $T= A_{47}$, and two years later, Xu {\em et al}~~\cite{XFWX} further showed that if $T=A_{47}$ and $\G$ is not normal, then $\G$ must be $5$-arc-transitive and up to isomorphism there are exactly two such graphs. Du {\em et al}~\cite{DFZh} showed that a connected pentavalent symmetric Cayley graph $\G$ over $T$ is either normal, or $\Aut(\G)$ contains a non-abelian simple normal subgroup $L$ such that $T\leq L$ and $(T,L)$ is one of $13$ possible pairs of non-abelian simple groups.

For $T$-symmetric graphs, Fang and Praeger~\cite{FangP,FP2} classified such graphs when $T$ is a Suzuki or Ree simple group acting transitively on the set of $2$-arcs of the graphs. For a connected cubic $T$-symmetric graph $\G$, it was proved by Li~\cite{CHLi} that either $\G$ is $T$-normal or $(T,\Aut(\G))=(A_7,A_8)$, $(A_7,S_8)$, $(A_7,2.A_8)$, $(A_{15},A_{16})$ or $(\GL(4,2),\rm AGL(4,2))$. Fang {\em et al}~\cite{FangLW} proved that none of the above five pairs can happen, that is, $T$ is always normal in $\Aut(\G)$. Du {\em et al}~\cite{DFZh} showed that a connected pentavalent $T$-symmetric graph $\G$ is either $T$-normal or $\Aut(\G)$ contains a non-abelian simple normal subgroup $L$ such that $T\leq L$ and $(T,L)$ is one of $17$ possible pairs of non-abelian simple groups.

Let $G$ be the characteristically simple group $T^l$ with $l\geq 1$. In this paper, we extend the above results on connected pentavalent $T$-vertex graphs to $G$-vertex graphs.

\begin{theorem}\label{theo=main}
Let $T$ be a non-abelian simple group and let $G=T^l$ with $l\geq 1$. Assume that every connected pentavalent symmetric $T$-vertex-transitive graph is $T$-normal. Then every connected pentavalent symmetric $G$-vertex-transitive graph is $G$-normal.
\end{theorem}

In 2011, Hua {\em et al}~\cite{huaF} proved that if every connected cubic symmetric $T$-vertex-transitive graph is $T$-normal, then every connected cubic symmetric $G$-vertex-transitive graph is $G$-normal. By Theorem~\ref{theo=main} and \cite[Theorem 1.1]{DFZh}, we have the following corollaries.

\begin{cor}\label{cor=vertex}
Let $T$ be a non-abelian simple group and let $G=T^l$ with $l\geq 1$. Then every connected pentavalent symmetric $G$-vertex-transitive graph is $G$-normal except for $T=\PSL(2,8)$, $\Omega_8^-(2)$ or $A_{n-1}$ with $n\geq 6$ and $n\di 2^9\cdot 3^2\cdot 5$.
\end{cor}

\begin{cor}\label{cor=arc}
Let $T$ be a non-abelian simple group and let $G=T^l$ with $l\geq 1$. Then every connected pentavalent $G$-symmetric graph is $G$-normal except for $T=A_{n-1}$ with $n=2\cdot 3$, $2^2\cdot 3$, $2^4$, $2^3\cdot 3$, $2^5$, $2^2\cdot 3^2$, $2^4\cdot 3$, $2^3\cdot 3^2$, $2^5\cdot 3$, $2^4\cdot 3^2$, $2^6\cdot 3$, $2^5\cdot 3^2$, $2^7\cdot 3$, $2^6\cdot 3^2$, $2^7\cdot 3^2$, $2^8\cdot 3^2$ or $2^9\cdot 3^2$.
\end{cor}

\begin{cor}\label{cor=regular}
Let $T$ be a non-abelian simple group and let $G=T^l$ with $l\geq 1$. Then every connected pentavalent symmetric $G$-regular graph is $G$-normal except for $T=\PSL(2,8)$, $\Omega_8^-(2)$ or $A_{n-1}$ with $n=2\cdot 3$, $2^3$, $3^2$, $2\cdot 5$, $2^2\cdot 3$, $2^2\cdot 5$, $2^3\cdot 3$, $2^3\cdot5$, $2\cdot3\cdot 5$, $2^4\cdot5$, $2^3\cdot 3\cdot5$, $2^4\cdot 3^2\cdot5$, $2^6\cdot 3\cdot5$, $2^5\cdot 3^2\cdot5$, $2^7\cdot 3\cdot5$, $2^6\cdot 3^2\cdot5$,  $2^7\cdot 3^2\cdot5$ or $2^9\cdot 3^2\cdot5$.
\end{cor}

\section{Preliminaries}

In this section, we describe some preliminary results which will be used later. The first one
is the vertex stabilizers of connected pentavalent symmetric graphs. By \cite[Theorem 1.1]{Guo}, we have the following proposition.

\begin{prop}\label{prop=stabilizer}
Let $\Gamma$ be a connected pentavalent $G$-symmetric graph with $v\in V(\Gamma)$. Then $G_v\cong \mathbb{Z}_5$, $D_{5}$, $D_{10}$, $F_{20}$, $F_{20}\times \mathbb{Z}_2$,
$ F_{20}\times \mathbb{Z}_4$, $A_5$, $S_5$, $A_4\times A_5$, $S_4\times S_5$, $(A_4\times A_5)\rtimes \mathbb{Z}_2$, $ {\rm ASL}(2,4)$, ${\rm AGL}(2,4)$,
${\rm A\Sigma L}(2,4)$,  ${\rm A\Gamma L}(2,4)$ or $\mathbb{Z}^6_2\rtimes {\rm \Gamma L}(2,4)$, where $F_{20}$ is the Frobenius group of order $20$,
 $A_4\rtimes \mathbb{Z}_2=S_4$ and $A_5\rtimes \mathbb{Z}_2=S_5$.
In particular, $|G_v|=5$, $2\cdot 5$, $2^2\cdot 5$, $2^2\cdot 5$, $2^3\cdot 5$, $2^4\cdot 5$, $2^2\cdot 3 \cdot 5$, $2^3\cdot 3 \cdot 5$, $2^4\cdot 3^2\cdot 5$, $2^6\cdot 3^2\cdot 5$, $2^5\cdot 3^2\cdot 5$, $2^6\cdot3 \cdot 5$, $2^6\cdot 3^2\cdot 5$, $2^7\cdot 3\cdot 5$, $2^7\cdot 3^2\cdot 5$ or $2^9\cdot 3^2 \cdot 5\cdot 5$, respectively.
\end{prop}

Connected pentavalent symmetric graphs admitting vertex-transitive non-abelian simple groups were classified in \cite{DFZh}.

\begin{prop}\label{prop=non-abelian}{\rm \cite[Theorem 1.1]{DFZh}}
Let $T$ be a non-abelian simple group and $\G$ a connected pentavalent symmetric $T$-vertex-transitive graph. Then either $T\unlhd \Aut(\G)$, or $T=\Omega^-_8(2), \PSL(2,8)$ or $A_{n-1}$ with $n\geq 6$ and $n\di 2^9\cdot3^2\cdot5$.
\end{prop}

The following is straightforward (also see the short proof of \cite[Lemma 3.2]{DFZh}).

\begin{prop}\label{prop=GH} Let $\G$ be a connected pentavalent symmetric $G$-vertex-transitive graph with $v\in V(\G)$ and let $A=\Aut(\G)$. If  $H\leq A$ and $GH\leq A$, then $|H|/|H\cap G|=|(GH)_v|/|G_v|\di 2^9\cdot 3^2 \cdot 5$, and if $\G$ is further $G$-symmetric then  $|H|/|H\cap G|\di 2^9\cdot 3^2$.
\end{prop}

The following proposition follows the classification of three-factor simple groups.

\begin{prop}\label{prop=235simplegroup}{\rm \cite[Theorem \MakeUppercase{\romannumeral1} ]{Huppert2}}
Let $G$ be a non-abelian simple $\{2,3,5\}$-group.
Then $G= A_5$, $A_6$ or $\rm PSU(4,2)$.
\end{prop}

By Guralnick~\cite[Theorem 1]{Guralnick}, we have the following proposition.

\begin{prop}\label{simplegroupwithp-powerindex}
Let $G$ be a non-abelian simple group with a subgroup $H$ such that $|G:H|=p^a$ with $p$ a prime and $a\geq1$. Then \begin{enumerate}
       \item [\rm(1)] $G=A_n$ and $H=A_{n-1}$ with $n=p^a$;
       \item [\rm(2)] $G=\PSL(2,11)$ and $H=A_5$ with $|G:H|=11$;
       \item [\rm(3)] $G=M_{23}$ and $H=M_{22}$ with $|G:H|=23$, or $G=M_{11}$ and $H=M_{10}$ with $|G:H|=11$;
       \item [\rm(4)] $G=\PSU(4,2)\cong \PSp(4,3)$ and $H$ is the parabolic subgroup of index $27$;
       \item [\rm(5)] $G=\PSL(n,q)$ and $H$ is the stabilizer of a line or hyperplane with $|G:H|=(q^n-1)/(q-1)=p^a$.
     \end{enumerate}
\end{prop}

By \cite[Theorem 1]{LiXu} and Proposition~\ref{simplegroupwithp-powerindex}, we have the following proposition.

\begin{prop}\label{prop=2a3b}
Let $G$ be a non-abelian simple group and $H$ a maximal subgroup of $G$ such that $|G:H|=2^a\cdot 3^b\geq 6$ with $0\leq a\leq 9$ and $0\leq b\leq 2$. Then $G$, $H$ and $|G:H|$ are listed in Table~\ref{table=1}.
\end{prop}

\begin{table}[ht]
\begin{center}
\begin{tabular}{|c|c|c||c|c|c|}

\hline
$G$              & $H$              & $|G:H|$              &$G$              & $H$                      & $|G:H|$           \\
\hline
$M_{11}$         & $\PSL(2,11)$     & $2^2\cdot3$          &$M_{12}$         & $M_{11}$                 & $2^2\cdot 3$         \\
\hline
$M_{24}$         & $M_{23}$         & $2^3\cdot 3$         &$\PSU(3,3)$      & $\PSL(2,7)$              & $2^2\cdot3^2$         \\
\hline
$A_9$            & $S_7$            & $2^2\cdot 3^2$       &$\PSU(4,2)$      & $S_6$                    & $2^2\cdot 3^2$         \\
\hline
$\PSp(6,2)$     & $S_8$            & $2^2\cdot 3^2$       & $M_{12}$        & $\PSL(2,11)$              & $2^4\cdot3^2$   \\
\hline
$\PSL(2,8)$     & $D_7$         & $2^2\cdot 3^2$       &$\PSL(3,3)$      & $\mz_{13}\rtimes \mz_3$   & $2^4\cdot 3^2$         \\
\hline
$\PSL(2,9)$     & $A_5$          & $2\cdot 3$       & $\PSL(2,p)$    & $\mz_{p}\rtimes \mz_{\frac{p-1}{2}}$     & $p+1=2^a\cdot 3^b$   \\
\hline
 $\PSL(2,8)$ & $\mz_2^3\rtimes\mz_7$  & $3^2$ & $A_n$      & $A_{n-1}$              & $n=2^a\cdot3^b$        \\
\hline
\end{tabular}
\end{center}
\vskip -0.5cm
\caption{{Non-abelian simple group pairs of index $2^a\cdot 3^b$}}\label{table=1}
\end{table}

Let $G$ be a group. The \emph{inner automorphism group} $\rm Inn(G)$ of $G$ is the group of automorphisms
of $G$ induced by conjugate action of elements in $G$, which is a normal subgroup in
the automorphism group $\Aut(G)$ of $G$. The quotient group $\Aut(G)/\rm Inn(G)$ is called
the \emph{outer automorphism group} of $G$. By the classification of finite simple groups, we have
the following proposition, which is the famous Schreier conjecture.

\begin{prop}\label{prop=auto out}{\rm \cite[Theorem 1.64]{Gorenstein}}
Every finite simple group has a solvable outer automorphism group.
\end{prop}

Baddeley and Praeger \cite{BP} considered almost simple groups containing a direct product
of at least two isomorphic non-abelian simple groups.
\begin{prop}\label{prop=T^k}{\rm \cite[Theorem 1.4]{BP}}
Let $H$ be an almost simple group, that is, $S\leq H\leq \Aut(S)$ for a non-abelian simple group $S$, and suppose that $H=AB$, where $A$ is a proper subgroup of $H$ not containing $S$, and $B\cong T^r$ for a non-abelian simple group $T$ and integer $k\geq2$. Then $S=A_n$ and $A\cap S = A_{n-1}$, where $n=|H: A|=|S:A\cap S| \geq 10$.
\end{prop}

Let $\G$ be a graph and $N \leq \Aut(\G)$. The \emph{quotient graph} $\G_N$ of $\G$ relative to $N$ is defined as the graph with vertices the orbits of $N$ on $V(\G)$ and with two orbits adjacent if there is an edge in $\G$ between these two orbits.

\begin{prop}\label{prop=atlesst3orbits}{\rm \cite[Theorem 9]{Lorimer}}
Let $\Gamma$ be a connected $G$-symmetric graph of prime valency, and let $N\unlhd G$ have at least three orbits on $V(\Gamma)$. Then $N$ is the kernel of $G$ on $V(\G_N)$, and semiregular on $V(\G)$. Furthermore, $\G_N$
is $G/N$-symmetric with $G/N \leq \Aut(\G_N)$.
\end{prop}

\section{Proof of Theorem~\ref{theo=main}}

Let $T$ be a non-abelian simple group and let $G=T^l$ with $l\geq 1$. Let $\G$ be a connected pentavalent symmetric $G$-vertex-transitive graph with $v\in V(\G)$ and let $A=\Aut(\G)$. We make the following assumption throughout this section.

\vskip 0.3cm

\f {\bf Assumption:} Every connected pentavalent symmetric $T$-vertex-transitive graph is $T$-normal.

Since the complete graph $K_6$ of order $6$ has automorphism group $S_6$ and is $A_5$-symmetric, we have $T\not=A_5$ by Assumption. Since $G$ is vertex-transitive and has no subgroup of index $2$, $\G$ is not bipartite.

To prove Theorem~\ref{theo=main}, we apply induction on $l$. It suffices to show that $G$ contains a minimal normal subgroup of $A$, and this is done in Lemmas~\ref{lem=minimal} and \ref{lem=insolvable} when $\rad(A)=1$, where $\rad(A)$ is the largest solvable normal subgroup of $A$. For $\rad(A)\not=1$, we need the fact $\rad(A)G=\rad(A)\times G$, which is proved in Lemmas~\ref{lem=RG} and \ref{lem=GP}.

\begin{lem}\label{lem=minimal}
Let $\G$ be $X$-symmetric with $G\leq X$, and let $X$ have a minimal normal subgroup that is a direct product of $T$. Then $G$ contains a minimal normal subgroup of $X$.
\end{lem}

\f {\bf Proof:} By Assumption, if $l=1$ then $G\unlhd A$ and the lemma is true. Assume $l\geq 2$.

Let $N$ be a minimal normal subgroup of $X$ such that $N=T^{m}$ for a positive integer $m$. Since $N\cap G\unlhd G$, we have $N\cap G=T^n$ with $n\leq m$ and $n\leq l$. Set $D=NG$. By Proposition~\ref{prop=GH}, $|N|/|N\cap G|=|T|^{m-n}=|D|/|G|=|D_v|/|G_v|\di2^9\cdot 3^2 \cdot 5$. Since $T$ is non-abelian simple, we have either $m-n=1$ or $m=n$.

Suppose $m-n=1$. Then $|T|\di2^9\cdot 3^2 \cdot 5$, and by Proposition~\ref{prop=235simplegroup}, $T=A_6$ as $T\not=A_5$. Thus, $|G|=|A_6|^l$, $|D|=|A_6|^{l+1}$ and $|D_v|/|G_v|=|A_6|$, implying $5\di |D_v|$ and $5\nmid |G_v|$. It follows that $\G$ is $D$-symmetric, but not $G$-symmetric.

Let $\bar{D}=D/G_D$, $\bar{G}=G/G_D$ and $\bar{N}=NG_D/G_D$, where $G_D$ is the largest normal subgroup of $D$ contained in $G$. Then $\bar{D}=\bar{G}\bar{N}$. Since $\G$ is $D$-symmetric, $D_v$ is primitive on the set of $5$ neighbours of $v$ in $\G$, and since $G_D\unlhd D$, we have either $(G_D)_v=1$ or $5\di |(G_D)_v|$. Thus,  $(G_D)_v=1$ because $\G$ is not $G$-symmetric. If $G=G_D$, then $G_v=1$ and $|D_v|=|T|=2^3\cdot 3^2\cdot5$, contrary to Proposition~\ref{prop=stabilizer}. Thus, $G/G_D=T^{l'}$ with $l'\geq 1$.

If $G_D$ is transitive on $V(\G)$, then $|G_v|=|G_v|/|(G_D)_v|=|G|/|G_D|=|T|^{l'}$. Since $T=A_6$, we have $5\di |G_v|$, a contradiction. If $G_D$ has two orbits on $V(\G)$, then $\G$ is bipartite, which is impossible.
Thus, $G_D$ has at least three orbits on $V(\G)$. By Proposition~\ref{prop=atlesst3orbits}, $\G_{G_D}$ is a connected pentavalent $G/G_D$-vertex-transitive and $D/G_D$-symmetric graph. If $l'=1$,
by Proposition \ref{prop=non-abelian}, $G/G_D\unlhd D/G_D$, that is, $G\unlhd D$, which is also impossible as otherwise $G_v=1$ and $|D_v|=|T|=2^3\cdot 3^2\cdot5$. Thus, $\bar{G}=G/G_D=T^{l'}$ with $l'\geq 2$ and $|\bar{D}|=|T|^{l'+1}$ as $|\bar{D}:\bar{G}|=|D:G|=|T|$.

Let $\bar{M}$ be a maximal subgroup of $\bar{D}$ containing $\bar{G}$. Then $\bar{G}\leq \bar{M}$ and $|T|^{l'}\leq |\bar{M}|<|T|^{l'+1}$. Consider the right multiplication of $\bar{D}$ on $[\bar{D}:\bar{M}]$, the set of right cosets of $\bar{M}$ in $\bar{D}$. Then $\bar{D}/\bar{M}_{\bar{D}}$ is a primitive permutation group on $[\bar{D}:\bar{M}]$ and $\bar{D}/\bar{M}_{\bar{D}}\leq S_{|T|}$, where $\bar{M}_{\bar{D}}$ is the largest normal subgroup of $\bar{D}$ contained in $\bar{M}$. Since
$|\bar{D}:\bar{M}|\di |T|=|A_6|$, the primitive group $\bar{D}/\bar{M}_{\bar{D}}$ has degree dividing $360$.

Let $\bar{L}$ be a minimal normal subgroup of $\bar{D}$.
Then $\bar{L}\cap \bar{N}=\bar{L}$ or $1$. For the former, $\bar{L}\leq \bar{N}$, and for the latter, $\bar{L}=\bar{L}/\bar{L}\cap \bar{N}\cong \bar{N}\bar{L}/\bar{N}\unlhd \bar{D}/\bar{N}\cong\bar{G}/\bar{G}\cap \bar{N}$.
In both cases, $\bar{L}$ is a direct product of $T$, and hence $\soc(\bar{D})=T^k$ for a positive integer $k$, where $\soc(\bar{D})$ is the product of all minimal normal subgroups of $\bar{D}$, that is, the socle of $\bar{D}$.

Assume $\bar{M}_{\bar{D}}\neq 1$. Then $\bar{D}$ has a minimal normal subgroup contained in $\bar{M}_{\bar{D}}$, say $\bar{K}$. Then  $\bar{K}\unlhd ~\soc(\bar{D})=T^k$. If $\bar{K}\nleq \bar{G}$, then $|T|^{l'+1}=|\bar{D}|>|\bar{M}|\geq |\bar{G}\bar{K}|=|\bar{G}|(|\bar{K}|/|\bar{K}\cap\bar{G}|) \geq |\bar{G}||T|=|T|^{l'+1}$, a contradiction. Thus, $\bar{K}\leq \bar{G}$, and we may assume $\bar{K}=K/G_D$, where $G_D\leq K\leq G$. Since $\bar{K}\unlhd \bar{D}$, we have $K\unlhd D$, and hence $K\leq G_D$. It follows that $K=G_D$ and $\bar{K}=1$, a contradiction.

Thus, $\bar{M}_{\bar{D}}= 1$ and so $\bar{D}$ is a primitive group on $[\bar{D}:\bar{M}]$ of degree dividing $|T|=360$. Since $\soc(\bar{D})=T^k=A_6^k$,  we have $k\leq 2$ by \cite[Section 7]{Colva}, and the primitivity implies that $|\bar{D}|\leq |\Aut(\soc(\bar{D}))|\leq 32|T|^2$, which is impossible because  $|\bar{D}|=|T|^{l'+1}\geq |T|^3$.

It follows that $m-n=0$ and hence $N\leq G$, as required.
\hfill\qed

\f {\bf Remark:} It is easy to check that Lemma \ref{lem=minimal} is also true if Assumption is replaced by the conditions that $l\geq 2$ and $T\neq A_5$.

\begin{lem}\label{lem=insolvable}
Let $G\leq X$ such that $\rm rad$$(X)=1$ and $\G$ is $X$-symmetric. Then $X$ has a minimal normal subgroup that is a direct product of $T$.
\end{lem}

\f {\bf Proof:} Suppose to the contrary that every minimal normal subgroup of $X$ is not a direct product of $T$. Then $G=T^l\ntrianglelefteq X$, and by Assumption, $l\geq 2$.

Let $L$ be the socle of $X$, the product of all minimal normal subgroups of $X$. Since $\rad(X)=1$, we have $L=S_1\times \cdots \times S_t$ with $t$ a positive integer, and $S_i\ncong T$
is a non-abelian simple group for each $1\leq i\leq t$. Furthermore, $L\cap G=T^r$ with $0\leq r\leq l$ as $L\cap G\unlhd G$. By Proposition~\ref{prop=GH}, $|L|/|L\cap G|=|(GL)_v|/|G_v|$
$\di 2^9\cdot 3^2 \cdot 5$.

For convenience, write $L\cap G=T^r=T_1\times\cdots \times T_r$. Let $P_i$ be the projection from $L\cap G$ to $S_i$ for a given $1\leq i\leq t$, that is, $P_i(x)=x_i$ for any $x=x_1\cdots x_t\in L\cap G$ with $x_i\in S_i$ ($1\leq i\leq t$). Then $P_i$ is a homomorphism from $L\cap G$ to $S_i$ and hence $(L\cap G)/K_i\cong I_i$, where $K_i$ is the kernel of $P_i$ and $I_i\leq S_i$ is the image of $L\cap G$ under $P_i$. Since $K_i\unlhd L\cap G$, we have $(L\cap G)/K_i\cong T^{r_i}=I_i$. Clearly, $K_1\leq S_2\times \cdots \times S_t$ and $K_1\cap \cdots\cap K_t=1$.
It follows $T^r=L\cap G=(L\cap G)/(K_1\cap\cdots\cap K_t)\lesssim (L\cap G)/K_1\times\cdots\times (L\cap G)/K_t\cong T^{r_1}\times \cdots\times T^{r_t}=I_1\times \cdots\times I_t$, and hence $(|S_1:I_1|\cdots |S_t:I_t|)\di |L:(L\cap G)|\di 2^9\cdot 3^2 \cdot 5$.

Assume $t\geq 5$. Then at least four $|S_i:I_i|$ are divisors of $2^9\cdot 3^2$, say $|S_i:I_i|=2^{a_i}\cdot 3^{b_i}$ with $2\leq i\leq 5$. Since $T\neq A_5$, Proposition \ref{prop=2a3b} implies $2^{a_i}\cdot 3^{b_i}\geq 2^3$ and $|S_2:I_2|\cdots|S_t:I_t|\geq 2^{12}$. If $5\nmid |S_1:I_1|$, then $|S_1:I_1|\geq 2^3$ and $|S_1:I_1|\cdots |S_t:I_t|\geq 2^{15}$.
Since $(|S_1:I_1|\cdots |S_t:I_t|)\di 2^9\cdot 3^2 \cdot 5$, we have $2^{15}\leq 2^9\cdot 3^2\cdot 5$, a contradiction. Thus, $5\di |S_1:I_1|$. If $|S_1:I_1|\not=5$ then $|S_1:I_1|\geq 2\cdot 5$ and hence $2^{13}\cdot 5\leq 2^9\cdot 3^2\cdot 5$, a contradiction. It follows $|S_1:I_1|=5$, which is also impossible by Proposition~\ref{simplegroupwithp-powerindex}. This yields $t\leq 4$.

Let $C=C_X(L)$ be the centralizer of $L$ in $X$. Then $C\cap L=1$ and hence $C=1$ because $L$ contains every minimal normal subgroup of $X$. Thus, $X=X/C\leq \Aut(L)$ and $G/G\cap L\cong GL/L\leq X/L\leq \Out(L)$. Since $t\leq 4$, Proposition~\ref{prop=auto out} implies that $\Out(L)$
is solvable, and hence $G/G\cap L=1$ and $G\leq L$. Since $G\ntrianglelefteq X$, we have $G<L$, that is, $G$ is a proper subgroup of $L$. It follows that $L_v\not=1$ and since $L\unlhd X$, we have $5\di |L_v|$ and so $L$ is symmetric.

If $t=1$, the Frattini argument implies that  $L=GL_v$, and by Proposition~\ref{prop=T^k}, $L_v=A_n$ with $n\geq 9$, which is impossible. Thus, $2\leq t\leq 4$. Since $(|S_1:I_1|\cdots |S_t:I_t|)\di 2^9\cdot 3^2\cdot 5$, there is at least one $|S_i:I_i|$ that is a divisor of $ 2^9\cdot 3^2$, say $|S_1:I_1|\di 2^9\cdot 3^2$. Since $I_1=T^{r_1}$,  Proposition~\ref{prop=2a3b} implies $I_1=T$, and since $I_1\cong G/K_1$, we have $K_1=T^{l-1}$. Set $\bar{S_1}=S_2\times \cdots \times S_t$. Then $L=S_1\times \bar{S_1}$ and $K_1=G\cap \bar{S_1}$ as $K_1\leq \bar{S_1}$ and $G=T^l\not= K_1$.

If $\bar{S_1}$ is transitive on $V(\G)$, then $|L_v|/|(\bar{S_1})_v|=|L|/|\bar{S_1}|=|S_1|\di 2^9\cdot 3^2\cdot5$ and so $S_1=A_5$ or $A_6$, forcing $I_1=T=A_5$, a contradiction.
If $\bar{S_1}$ has two orbits, then $\G$ is bipartite, which is also impossible. Thus, $\bar{S_1}$ has at least three orbits on $V(\G)$ and Proposition~\ref{prop=atlesst3orbits}
implies that $\G_{\bar{S_1}}$ is a connected pentavalent $G\bar{S_1}/\bar{S_1}$-vertex-transitive
and $L/\bar{S_1}$-symmetric graph. Since $G\bar{S_1}/\bar{S_1}\cong G/G\cap \bar{S_1}=T$,
by Assumption, $T\cong G\bar{S_1}/\bar{S_1}\unlhd L/\bar{S_1}\cong S_1$, which is impossible because
$S_1$ is a non-abelian simple group. This completes the proof.
\hfill\qed

\f {\bf Remark:} Let
 $$\begin{array}{cl}
             \Delta= & \{\PSL(2,8), \ \Omega_8^-(2),\ A_{n-1}\ |\ n\geq 6,\ n\di 2^9\cdot 3^2\cdot 5\}\\
              \Delta_1= & \{\PSL(2,8),\ \Omega_8^-(2),\ A_{n-1}\ |\ n=2\cdot 3, 2^3, 3^2, 2^2\cdot 3, 2\cdot 5, 2^2\cdot 5, 2^3\cdot 3, \\
              &  2^3\cdot5, 2\cdot3\cdot 5, 2^5\cdot 3^2\cdot5, 2^7\cdot 3\cdot5, 2^6\cdot 3^2\cdot5, 2^7\cdot 3^2\cdot5, 2^9\cdot 3^2\cdot5\}.
\end{array}$$
Then Lemma~\ref{lem=insolvable} is true if Assumption is replaced by the conditions that $l\geq2$ and $T=A_{n-1}\in \Delta-\Delta_1$.
To prove it, we only need to change the last paragraph in the proof of Lemma~\ref{lem=insolvable} as following.

Since $|S_1:I_1|\di 2^9\cdot3^2$, Proposition \ref{prop=2a3b} implies that $I_1=T$, and since $T\in \Delta-\Delta_1$, we have $(S_1,I_1)=(A_n,A_{n-1})$ with $n=2^i$, $2^j\cdot3$ or $2^k\cdot 3^2$, where $4\leq i,j\leq 9$ and $1\leq k\leq 9$.
Recall that $(|S_1:I_1|\cdots |S_t:I_t|)\di 2^9\cdot 3^2\cdot 5$ and $2\leq t\leq 4$. Then for each $2\leq i\leq t$, we have $|S_i:I_i|\di 2^5\cdot 3^2\cdot5$ or $2^8\cdot5$
and hence $S_i$ is a primitive permutation group of degree dividing $1440$ or $1280$. By \cite[Section 7]{Colva}, $(S_i:I_i)=(A_n,A_{n-1})$, which implies that $t=2$ and $(S_2,S_1,T)=(A_{n},A_n, A_{n-1})$ with $n=2^4$ or $2^4\cdot3$. Thus, $G=T^2<L=S_1\times S_2$, $G\cap S_2=T$ (note that $\bar{S_1}=S_2$), $|L_v|/|G_v|=|L|/|G|=n^2$ and $2^8\di |L_v|$. By Proposition~\ref{prop=stabilizer}, $|L_v|=2^9\cdot3^2\cdot5$ and so $|G_v|=2\cdot3^2\cdot5$ or $2\cdot5$. The former is impossible by Proposition~\ref{prop=stabilizer}. For the latter, $|(GS_2)_v|/|G_v|=|GS_2|/|G|=|S_2|/|S_2\cap G|=|S_2|/|T|=2^4\cdot3$ and $|(GS_2)_v|=2^5\cdot 3\cdot5$, which is also impossible by Proposition~\ref{prop=stabilizer}.

\begin{lem}\label{lem=RG}
Assume $5\di \rm |rad$$(A)|$. Then $\rad(A)G=\rad(A)\times G$.
\end{lem}

\f {\bf Proof:} Since $\rad(A)$ is solvable, $\rad(A)\cap G=1$. If $l=1$ then by Assumption, $G\unlhd A$, and hence $\rad(A)G=\rad(A)\times G$. Thus, we may assume $l\geq 2$.

Set $B=\rad(A)G$. Then $|B|=|\rad(A)||G|$, and by Proposition~\ref{prop=GH}, $|\rad(A)|=|B|/|G|=|B_v|/|G_v|\di 2^9\cdot 3^2 \cdot 5$. Since $5 \di |\rad(A)|$, we have $5\di |B_v|$ and $5\nmid |G_v|$, that is, $\G$ is $B$-symmetric, but not $G$-symmetric. In particular, $G_v$ is a $\{2,3\}$-group and $5^2\nmid |\rad(A)|$.

Since $\rad(A)$ is a solvable $\{2,3,5\}$-group, $\rad(A)$ has a Hall $\{2,3\}$-subgroup, say $H$.
Set $\Omega=\{H^r\ |\ r\in \rad(A)\}$. By~\cite{Hall}, all Hall $\{2,3\}$-subgroups of
$\rad(A)$ are conjugate and so the conjugate action of $B$ on $\Omega$ is transitive. Let $K$ be the kernel of the action of $G$ on $\Omega$. Since $5^2\nmid \rad(A)$, we have $|\Omega|=1$ or $5$, and hence $G/K=T^r\leq S_5$. Since $T\not=A_5$, we have $G=K$ and so $G$ fixes $H$. It follows $H\unlhd GH\leq B$ and $G\cap H=1$.

Set $Y=GH$ and $\Delta=\{Yb\ |\ b \in B\}$. Then $|\Delta|=|B:Y|=5$. Let $Y_B$ be the kernel of the right multiplication action of $B$ on $\Delta$. Then $Y_B$ is the largest normal subgroup of $B$ contained in $Y$. It follows $B/Y_B \leq S_5$. Suppose $G\nleq Y_B$. Then $G/G\cap Y_B\cong GY_B/Y_B\leq B/Y_B\leq S_5$, and so $G/G\cap Y_B=T=A_5$, a contradiction. Thus, $G\leq Y_B$.

Since $|Y_v|/|G_v|=|Y|/|G|=|H|$ and $G_v$ is a $\{2,3\}$-group, $Y_v$ is a $\{2,3\}$-group. Since $\G$ is
$B$-symmetric, $B_v$ is primitive on the set of the $5$ neighbors of $v$ in $\G$, and since $Y_B\unlhd B$, we have either $(Y_B)_v=1$ or $5\di |(Y_B)_v|$. The latter cannot happen as $Y_v$ is a $\{2,3\}$-group. Thus, $(Y_B)_v=1$, and the Frattini argument implies $Y_B=G(Y_B)_v=G$. Thus, $G\unlhd B$ and $B=\rad(A)\times G$.
\hfill\qed

\begin{lem}\label{lem=GP} Let $G\leq X\leq A$ and let $\G$ be $G$-symmetric. Then $\rad(X)G=\rad(X)\times G$.
\end{lem}

\f {\bf Proof:} The lemma is true for $|\rad(X)|=1$. In what follows we assume that $|\rad(X)|\neq 1$. Recall that $G=T^l$ with $l\geq 1$. Since $\G$ is $G$-symmetric, $5\di |G_v|$ and hence $5\di |T|$.
Set $B=\rad(X)G$. Then $G\cap \rad(X)=1$ as $\rad(X)$ is solvable, and so $|B|=|\rad(X)||G|$. By Proposition~\ref{prop=GH}, $|\rad(X)|=|B_v|/|G_v|\di 2^9\cdot 3^2 $.

Let $N$ be a minimal normal subgroup of $X$ contained
in $\rad(X)$. Then $N\cong \mz_2^s$ or $N=\mz_3^t$ for some $1\leq s\leq 9$ or $1\leq t\leq 2$, and hence $\Aut(N)\leq \GL(s,2)$ or $\GL(t,3)$.

First we claim $D:=GN=G \times N$. Consider the conjugate action of
$G$ on $N$ and let $K$ be the kernel of $G$ in this action. Then $K=T^r\unlhd G$ for some $r\leq l$, and $G/K\leq \Aut(N)$. It is easy to see that $D=G\times N$ if and only if $G=K$.

Suppose $G\neq K$. Then $r<l$ and $G/K=T^{l-r}$ is insolvable. In particular, $\Aut(N)$ is insolvable. Recall that $\Aut(N)\leq \GL(s,2)$ or $\GL(t,3)$. It follows that $N\cong \mz_2^s$ with $3\leq s\leq 9$ as both $\GL(2,2)$ and $\GL(2,3)$ are solvable.

Note that $G\cap N=1$. By Proposition~\ref{prop=GH}, $|N|=|D|/|G|=|D_v|/|G_v|$, and so $|D_v|/|G_v|$ is a $2$-power. If $G_v\cong \mz_5$ then Proposition~\ref{prop=stabilizer} implies that $D_v\cong D_5$, $D_{10}$, $F_{20}$, $F_{20}\times \mz_2$ or $F_{20}\times \mz_4$. It follows that $|N|\di 2^5$, and this is always true by checking all other possible cases for $G_v$ in Proposition~\ref{prop=stabilizer}. This means that $N\cong \mz_2^s$ with $3\leq s\leq 5$ and hence $G/K\leq \Aut(N)\leq \PSL(5,2)$. Since   $5^2\nmid |\PSL(5,2)|$, we have $l-r=1$ and $G/K=T$.

If $K$ has at least three orbits, Proposition~\ref{prop=atlesst3orbits} implies that $\G_K$ is a connected pentavalent $G/K$-symmetric graph, and by Assumption, $G/K\unlhd D/K$. It follows that $G\unlhd D$ and hence $G=K$, a contradiction. Thus, $K$ has one or two orbits. If $K$ has two orbits, then $\G$ is bipartite, a contradiction. This yields that $K$ is transitive. Then $|G_v|/|K_v|=|G|/|K|=|T|$, and since $K\unlhd G$ and $\G$ is $G$-symmetric, we have either $K_v=1$ or $5\di |K_v|$. If $5\di |K_v|$, then $5\nmid |T|$ because $|G_v|/|K_v|=|T|$, a contradiction. It follows $K_v=1$. Let $L=KN$. Then $L\unlhd D$ and so $L_v=1$ or $5\di |L_v|$. On the other hand, $|L_v|=|L_v|/|K_v|=|L|/|K|=|N|=2^s$ ($3\leq s\leq 5$), a contradiction. Therefore, $D=G \times N$, as claimed.

Now we finish the proof by induction on $\rm |rad(X)|$. Since $5\nmid |N|$,  we have $N_v=1$. Assume that $N$ has one or two orbits on $V(\G)$. Then $5\nmid |V(\G)|$, and so $5^2\nmid |G|$, which implies that $G$ is a simple group. By Assumption, $G\unlhd X$ and hence $\rad(X)G=G\times \rad(X)$. Assume that $N$ has at least three orbits. By Proposition~\ref{prop=atlesst3orbits}, $\G_N$ is a connected pentavalent $GN/N$-symmetric graph. Note that $GN/N\cong T^l$ and $GN/N\leq X/N$. Since $\rad(X/N)=\rad(X)/N$, we have $|\rad(X/N)|<|\rad(X)|$ and the  inductive hypothesis implies that $\rad(X)G/N=\rad(X/N)\cdot GN/N=\rad(X)/N\times GN/N$. Thus, $GN\unlhd \rad(X)G$. Since $GN=G\times N$, $G$ is characteristic in $GN$ and hence $G$ is normal in $\rad(X)G$. It follows $\rad(X)G=\rad(X)\times G$.
\hfill\qed

\medskip
Now, we are ready to prove Theorem \ref{theo=main}.
\medskip

\f {\bf The proof of Theorem \ref{theo=main}:} Recall that $G=T^l$ and $\G$ is a connected pentavalent symmetric $G$-vertex-transitive graph with  $A=\Aut(\G)$. We apply induction on $l$. If $l=1$, Theorem \ref{theo=main} is true by Assumption. Assume $l\geq 2$. Let $R$ be the radical of $A$, and set $B=RG$. Then $G\cap R=1$ and $|R|=|B|/|G|=|B_v|/|G_v|$.

If $R$ is transitive on $V(\G)$, Proposition~\ref{prop=GH} implies $|G|=|GR|/|R|=|(GR)_v|/|R_v|\di 2^9\cdot 3^2\cdot 5$, which is impossible because $l\geq 2$. Furthermore, $R$ cannot have two orbits as $\G$ is not bipartite. Thus, $R$ has at least three orbits. By Proposition \ref{prop=atlesst3orbits}, $\G_R$ is a connected pentavalent $B/R\cong G=T^l$-vertex-transitive and $A/R$-symmetric graph.
Since $R$ is the largest solvable normal subgroup of $A$, we have $\rm rad$$(A/R)=1$.
By Lemmas~\ref{lem=insolvable} and \ref{lem=minimal}, $A/R$ has a minimal normal subgroup $M/R=T^r$ contained in $B/R$ with $1\leq r\leq l$. It follows that $R\leq M\leq B$ and $M\unlhd A$. In particular, $R\not=M$.

Assume $B=R\times G$. Then $M=M\cap RG=R(M\cap G)=R\times (M\cap G)$ and $M\cap G$ is characteristic in $M$.
It follows $M\cap G\unlhd A$ and $|M\cap G|=|T|^r$. If $M\cap G$ is transitive, then $|G_v|/|(M\cap G)_v|=|G|/|M\cap G|=|T|^{l-r}\di 2^9\cdot3^2\cdot5$. It follows $l-r=0$ or $l-r=1$. The former implies $G=M\cap G\unlhd A$ and we are done. For the latter, $T=A_6$ as $T\neq A_5$ and hence $5\di |G_v|$ and $5\nmid |(M\cap G)_v|$, that is,
$\G$ is $G$-symmetric, but not $(M\cap G)$-symmetric. Since $M\cap G\unlhd G$, we have $(M\cap G)_v=1$ and $|G_v|=|T|=2^3\cdot 3^2\cdot5$, which is impossible by Proposition~\ref{prop=stabilizer}. If $M\cap G$ has two orbits on $V(\G)$, then $\G$ is bipartite, which is also impossible. Thus, $M\cap G$ has at least three orbits on $V(\G)$ and Proposition \ref{prop=atlesst3orbits} implies that $\G_{M\cap G}$ is a connected pentavalent $G/(M\cap G)$-vertex-transitive and $A/(M\cap G)$-symmetric graph. By inductive hypothesis, $G/(M\cap G)\unlhd A/(M\cap G)$ and so $G\unlhd A$.

By the above paragraph, to finish the proof we only need to show $B=R\times G$, and to do this, we now claim $B\unlhd A$. If $M$ is transitive, then $|B_v|/|M_v|=|B|/|M|=|B/R|/|M/R|=|T|^{l-r}\di2^9\cdot 3^2\cdot5$ and so $l-r=0$ or $l-r=1$. For the former, $B=M\unlhd A$, as claimed. For the latter, $T=A_6$ as $T\not=A_5$, and so $|B_v|/|M_v|=|A_6|$. It follows $M_v=1$ and $|B_v|=|T|=2^3\cdot 3^2\cdot5$, which is impossible by Proposition~\ref{prop=stabilizer}. If $M$ has two orbits, then $\G$ is bipartite, which is also impossible. Thus, $M$ has at least three orbits, and by Proposition~\ref{prop=atlesst3orbits}, $\G_M$ is a connected pentavalent $B/M$-vertex-transitive and $A/M$-symmetric graph. Since $B/M\cong (B/R)/(M/R)=T^{l-r}$ with $l-r<l$, the inductive hypothesis implies that $B/M\unlhd A/M$ and so $B\unlhd A$, as claimed.

Since $G\leq B\unlhd A$ and $A$ is symmetric, we have $B_v=1$ or $5\di |B_v|$. If $B_v=1$, the Frattini argument implies $B=GB_v=G$, and hence $R=1$ and $B=R\times G$. If $5\di |B_v|$, then $5\di |R|$ or $5\di |G_v|$ as $|B_v|=|G_v||R|$.
By Lemmas \ref{lem=RG} and \ref{lem=GP}, $B=$ $R\times G$. \hfill\qed

\f {\bf The proof of Corollary~\ref{cor=vertex}:}
This follows from Theorem~\ref{theo=main} and \cite[Theorem 1.1]{DFZh}. \hfill\qed

\f {\bf The proof of Corollary \ref{cor=arc}:} It is easy to check that Lemmas~\ref{lem=minimal}-\ref{lem=GP} are true if Assumption is replaced by the following assumption:

\medskip
\f {\bf Assumption 1:} Every connected pentavalent $T$-symmetric graph is $T$-normal.

\medskip
Then a similar proof to the proof of Theorem \ref{theo=main} implies that
under Assumption~1, every connected pentavalent $G$-symmetric graph is $G$-normal, and by \cite[Corollary 1.2]{DFZh}, we have  Corollary~\ref{cor=arc}. \hfill\qed

\f {\bf The proof of Corollary \ref{cor=regular}:} By assumption of  Corollary~\ref{cor=regular}, $\G$ is a connected pentavalent symmetric $G$-regular graph with $G=T^l$, and hence it is  $G$-vertex-transitive. Recall that
$$\begin{array}{cl}
             \Delta= & \{\PSL(2,8), \ \Omega_8^-(2),\ A_{n-1}\ |\ n\geq 6, \ n\di 2^9\cdot 3^2\cdot 5\},\\
              \Delta_1= & \{\PSL(2,8),\ \Omega_8^-(2),\ A_{n-1}\ |\ n=2\cdot 3, 2^3, 3^2, 2^2\cdot 3, 2\cdot 5, 2^2\cdot 5, 2^3\cdot 3, \\
              &  2^3\cdot5, 2\cdot3\cdot 5, 2^5\cdot 3^2\cdot5, 2^7\cdot 3\cdot5, 2^6\cdot 3^2\cdot5, 2^7\cdot 3^2\cdot5, 2^9\cdot 3^2\cdot5\}.
\end{array}$$
By Corollary~\ref{cor=vertex}, $\G$ is $G$-normal except for $T\in \Delta$, and to prove Corollary~\ref{cor=regular}, it suffices to show that $\G$ is $G$-normal for $T\in \Delta-\Delta_1$. In what follows we always assume $T\in \Delta-\Delta_1$.

Clearly, $T\not=A_5$. If $l=1$, Corollary~\ref{cor=regular} is true by \cite[Corollary 1.3]{DFZh}. Assume $l\geq 2$. Then Lemmas~\ref{lem=minimal} and \ref{lem=insolvable} are true by the Remarks of these two lemmas. Then a similar argument to the proof of Theorem~\ref{theo=main} (the first three paragraphs) implies that $\G$ is $G$-normal if $RG=R\times G$, where $R$ is the radical of $A=\Aut(\G)$.
Set $B=RG$. To finish the proof, we only need to show $B=R\times G$.

Since $R$ is solvable, $G\cap R=1$ and $|B|=|G||R|$, and by Proposition~\ref{prop=GH}, $|R|=|R|/|R\cap G|\di 2^9\cdot 3^2\cdot 5$. We may write $|R|=2^m\cdot3^h\cdot 5^k$, where $0\leq m\leq 9$,
$0\leq h\leq 2$ and $0\leq k\leq 1$. Since $R$ is solvable, there exists a series of principle subgroups of $B$:
\begin{center}
$B>R=R_{s}>\cdots R_1>R_0=1$
\end{center}
such that $R_i\unlhd B$ and $R_{i+1}/R_i$ is an elementary abelian $r$-group with $0\leq i\leq s-1$, where  $r=2$, $3$ or $5$. Clearly, $G$ has a natural action on $B_{i+1}/B_i$ by conjugation.

Suppose to the contrary that $B\neq R\times G$. Then there exists some $0\leq j\leq s-1$ such that $GR_j=G\times R_j$, but $GR_{j+1}\not=G\times R_{j+1}$. If $G$ acts trivially on $R_{j+1}/R_j$ by conjugation, then $[GR_j/R_j,R_{j+1}/R_j]=1$. Furthermore, since $GR_j/R_j\cong G=T^l$, we have $(GR_j/R_j)\cap (R_{j+1}/R_j)=1$, and since $|GR_{j+1}/R_j|=|GR_{j+1}/R_{j+1}||R_{j+1}/R_j|=
|G||R_{j+1}/R_j|=|GR_j/R_j||R_{j+1}/R_j|$, we have $GR_{j+1}/R_j=GR_j/R_j\times R_{j+1}/R_j$. In particular, $GR_j\unlhd GR_{j+1}$ and so $G\unlhd GR_{j+1}$ because $GR_j=G\times R_j$ implies that $G$ is characteristic in $GR_j$. It follows that $GR_{j+1}=G\times R_{j+1}$, a contradiction. Thus, $G$ acts non-trivially on $R_{j+1}/R_j$. Let $K$ be the kernel of this action. Then $G/K=T^{l'}$ with $l'\geq 1$.

Since $|R|=2^m\cdot3^h\cdot 5^k$ with $0\leq m\leq 9$,
$0\leq h\leq 2$ and $0\leq k\leq 1$, $R_{j+1}/R_j$ is an elementary abelian group of order $r^\ell$, where $\ell\leq m\leq 9$ for $r=2$, $\ell\leq 2$ for $r=3$ and $\ell\leq 1$ for $r=5$. Then $G/K\leq \GL(\ell,r)$. Since $\GL(2,3)$ and $\GL(1,5)$ are solvable, we have
$G/K\leq \GL(9,2)$, and since $T\in \Delta-\Delta_1$, we have $11\di |T|$ and hence $11\di |G/K|$. It follows $11\di |\GL(9,2)|$, a contraction. This implies $B=R\times G$, as required.
\hfill\qed

\medskip
\f {\bf Acknowledgements:} This work was supported by the National Natural Science Foundation of China (11571035) and by the 111 Project of China (B16002).


\begin{thebibliography}{99}
\bibitem{BP}
R. W. Baddeley and C. E. Praeger, On primitive overgroups of quasiprimitive permutation
groups, J. Alg. 263 (2003) 294-344.
\bibitem{DFZh}
J.-L. Du, Y.-Q. Feng, J.-X. Zhou, Pentavalent symmetric graphs admitting vertex-transitive non-abelian simple groups, Europ. J. Combin. 63 (2017) 132-145.
\bibitem{FangLW}
X.G. Fang, L.J. Jia, J. Wang, On the automorphism groups of symmetric graphs admitting an almost simple group, Europ. J. Combin. 29 (2008) 1467-1472.
\bibitem{FangP}
X.G. Fang, C.E. Praeger, Finite two-arc transitive graphs admitting a Suzuki simple group,
J. Algebraic Combin. 27 (1999) 3727-3754.
\bibitem{FP2}
X.G. Fang and C.E. Praeger, Finite two-arc-transitive graphs admitting a Ree simple group, Comm. Algebra 27 (1999) 3755-3769.
\bibitem{Gorenstein}
D. Gorenstein, Finite simple groups (Plenum Press, New York, 1982).
\bibitem{Guo}
S.-T. Guo, Y.-Q. Feng, A note on pentavalent $s$-transitive graphs,
Discrete Math. 312 (2012) 2214-2216.
\bibitem{Guralnick}
R.M. Guralnick, Subgroups of prime power index in a simple group,
J Algebra. 81 (1983) 304-311.
\bibitem{Hall}
P. Hall, A note on soluble groups, J. London Math. Soc. 3 (1928) 98-105.
\bibitem{huaF}
X.-H. Hua, Y.-Q. Feng, Cubic graphs admitting transitive non-abelian characterically simple groups,
Proceedings of the Edinburgh Mathematical Soc. 54 (2011) 113-123.
\bibitem{Huppert2}
B. Huppert, W. Lempken, Simple groups of order divisible by at most four primes,
Proc. of the F. Scorina Gemel State University 16 (2000) 64-75.
\bibitem{CHLi}
C.H. Li, Isomorphisms of finite Cayley graphs,
Ph.D. Thesis, The University of Western Australia, 1996.
\bibitem{LiXu}
X.H. Li, M.Y. Xu, The primitive permutation groups of degree $2^a\cdot3^b$,
Arch Math, 86 (2006) 385-391.
\bibitem{Lorimer}
P. Lorimer, Vertex-transitive graphs: symmetric graphs of prime valency,
J. Graph Theory 8 (1984) 55-68.
\bibitem{Colva}
C. M. Roney-Dougal, The primitive permutation groups of degree less than 2500, J.
Alg. 292 (2005) 154-183.
\bibitem{XFWX2005}
S.J. Xu, X.G. Fang, J. Wang, M.Y. Xu, On cubic $s$-arc transitive Cayley graphs of finite simple groups,
Europ. J. Combin. 26 (2005) 133-143.
\bibitem{XFWX}
S.J. Xu, X.G. Fang, J. Wang, M.Y. Xu, $5$-arc transitive cubic Cayley graphs on finite simple groups,
Europ. J. Combin. 28 (2007) 1023-1036.
\end{thebibliography}
\end{document}